\documentclass[nopreprintline,12pt]{elsarticle}

\usepackage{amsthm}
\usepackage{amsmath}
\usepackage{amssymb}
\usepackage{amsfonts}
\usepackage{dsfont}
\usepackage{mathrsfs}

\usepackage{color}
\usepackage{xcolor}

\usepackage[colorlinks=true,%
  linkcolor=black,anchorcolor=black,%
  citecolor=black,filecolor=black,%
  menucolor=black,%
  urlcolor=blue]{hyperref} 

\usepackage{ulem}

\usepackage{graphicx,subfigure,psfrag}
\usepackage{tikz}
\usepackage{wrapfig} 
\usepackage{multicol,colordvi,epsfig,float}
\usepackage{multirow}
\usepackage{enumerate}
\usepackage{indentfirst}
\usepackage{pgfplots}
\pgfplotsset{width=9cm,compat=1.15}

\newcommand{\D}[1]{\displaystyle{#1}}

\usepackage{verbatim}

\newtheorem{teo}{Theorem}

\newtheorem{defi}{Definition}

\theoremstyle{remark}

\newtheorem{obsv}{\bf Remark}

\newtheorem{proc}{\bf Procedure}

\journal{arXiv}

\begin{document}

\begin{frontmatter}



\title{A graph-based approach to customer segmentation using the RFM model}

 \author[label1]{André Luiz Corrêa Vianna Filho}
 \ead{andrevf@ufpr.br}
 \author[label1,label2]{Leonardo de Lima}
 \ead{leonardo.delima@ufpr.br}
 \author[label1,label3]{Mariana Kleina\corref{cor1}}
 \ead{marianakleina@ufpr.br}
 \cortext[cor1]{Corresponding author}
 \affiliation[label1]{organization={Graduate Program in Production Engineering, Federal University of Parana}}
\affiliation[label2]{organization={Department of Business Administration, Federal University of Parana}}
\affiliation[label3]{organization={Department of Production Engineering, Federal University of Parana}}
%%  \affiliation[label3]{organization={Federal Institute of Parana},
%%             addressline={Engenheiro Tourinho Street},
 %%            city={Campo Largo},
%%             postcode={83607-140},
%%             state={Parana},
%%             country={Brazil}}

%%
%%  \affiliation[label3]{organization={Administration School of the Federal University of Paraná},
%             addressline={Prefeito Lothário Meissner Avenue},
%             city={Curitiba},
%             postcode={80210-170},
%             state={Paraná},
%             country={Brazil}}
%\author{} %% Author name

%%\affiliation{organization={Department of Production Engineering of the Federal University of Paraná (UFPR)},%Department and Organization
            %addressline={Francisco Heráclito dos Santos Street}, 
            %city={Curitiba},
            %postcode={81530-000}, 
            %state={Paraná},
            %country={Brazil}}

%% Abstract
\begin{abstract}
  The present article proposes a graph-based approach to customer segmentation, combining the RFM analysis with the classical optimization max-$k$-cut problem. We consider each customer as a vertex of a weighted graph, and the edge weights are given by the distances between the vectors corresponding to the $(R,F,M)$-scores of the customers. We design a procedure to build a reduced graph with fewer vertices and edges, and the customer segmentation is obtained by solving the max-$k$-cut for this reduced graph. We prove that the optimal objective function values of the original and the reduced problems are equal. Additionally, we show that an optimal solution to the original problem can be easily obtained from an optimal solution to the reduced problem, which provides an advantage in dealing with computational complexity in large instances. Applying our methodology to a real customer dataset allowed us to identify distinct behaviors between groups and analyze their meaning and value from a business perspective.

\end{abstract}



  %\item A graph-based approach is proposed to solve the customer segmentation problem.
  %\item It combines the RFM analysis and the max-$k$-cut problem
  %\item We propose a graph-based clustering method that combines the RFM analysis and the max-$k$-cut problem and apply it to the customer segmentation problem.
  %\item A reduced problem, which provides the same optimal objective function value of the max-$k$-cut formulation to the original problem, makes the size of the graph to be partitioned by the max-$k$-cut problem independent of the number of customers to be segmented.
  %\item Moreover, we prove that a solution to the original problem can be easily obtained from a solution to the reduced problem.
  %\item A computational experiment was conducted on a large dataset, and the results are analyzed from a business perspective.

%% Keywords

\end{frontmatter}




%%%%%%%%%%%%%%%%%%%%%%%%%%%%%%%%%%%%%%%%%%%%%%%%%%%%%%%%%%%%%%%%%%%%%%%%%%%%%%%%%%%%%
%%%%%%%%%%%%%%%%%%%%%%%%%%%%%%%%%%%%%%%%%%%%%%%%%%%%%%%%%%%%%%%%%%%%%%%%%%%%%%%%%%%%%
%%%%%%%%%%%%%%%%%%%%%%%%%%%%%%%%%%%%%%%%%%%%%%%%%%%%%%%%%%%%%%%%%%%%%%%%%%%%%%%%%%%%%
%%%%%%%%%%%%%%%%%%%%%%%%%%%%%%%%%%%%%%%%%%%%%%%%%%%%%%%%%%%%%%%%%%%%%%%%%%%%%%%%%%%%%
%%%%%%%%%%%%%%%%%%%%%%%%%%%%%%%%%%%%%%%%%%%%%%%%%%%%%%%%%%%%%%%%%%%%%%%%%%%%%%%%%%%%%

\section{Introduction}

Customers are the primary source of profit for a company, and the satisfaction of their requirements is one of the key factors for the success of a business \cite{cheng2009classifying}. In this sense, Customer Relationship Management (CRM), described as a set of actions taken by an enterprise to understand and influence customer behavior to enhance customer acquisition, customer retention, customer loyalty and customer profitability \cite{swift2001accelerating}, plays a central role in a company's activities, helping it improve the loyalty of the existing customers and attract new ones.

According to \cite{ngai2009application}, and references cited therein, CRM consists of four dimensions: (i) customer identification, (ii) customer attraction, (iii) customer retention, and (iv) customer development. Our method is a tool to assist the \textit{customer identification} phase. This first CRM phase involves identifying groups of customers, including, for instance, the most profitable ones and those that are being lost to competitors. It is crucial because understanding their customers helps businesses develop targeted marketing strategies and endure in competitive environments \cite{rungruang2024rfm}. To support this phase, companies collect and store a great amount of data about their customers, and the challenge is then to extract, from this raw data, the useful information that will assist decision-making. An often-used tool for \textit{customer identification} is \textit{customer segmentation}, which counts with several methods, including clustering \cite{rungruang2024rfm}. 

In the present article, we propose a graph-based clustering method to be applied to \textit{customer segmentation}. To segment customers into clusters, we need to identify characteristics that allow us to distinguish them. To this end, we use the \textit{customer value analysis} method known as RFM analysis \cite{cheng2009classifying,hughes1994strategic}. Although there is a wealth of information about customers that companies can collect, the advantage of the RFM analysis, according to \cite{cheng2009classifying}, is that it identifies the characteristics of customers and differentiates them using only three variables, described in the sequel \cite{cheng2009classifying,christy2021rfm}: 
\begin{enumerate}
  \item Recency ($R_{i}$): the time interval between the latest consuming behavior of customer $i$ and the present. The shorter this interval is, the bigger the value of $R_{i}$.
  \item Frequency ($F_{i}$): the number of customer $i$ purchases in a fixed period. The higher this number is, the higher the value of $F_{i}$.
  \item Monetary ($M_{i}$): the total amount of money spent by customer $i$ over the whole time period considered. The higher this value is, the higher the value of $M_i$.
\end{enumerate}

Each variable is ordered in ascending order and divided into $T$ equal groups. Then, for each customer $i$, we assign an integer value ranging from $1$ to $T$ for each variable $R_{i}$, $F_{i}$ and $M_{i}$. So, each customer $i$ is associated with a three-dimensional score $(R_i, F_i, M_i)$ that we refer to as $(R,F,M)$-score of customer $i$. Note that there are at most $T^3$ distinct $(R,F,M)$-scores.

The customer segmentation problem consists of partitioning the dataset of an organization's customers into groups so that customers of the same group share common characteristics. A widely used strategy to solve this problem is to combine the RFM analysis with a clustering method \cite{cheng2009classifying,chen2012data,christy2021rfm,dougan2018customer,rungruang2024rfm}. There are several heuristic algorithms for clustering, and $K$-means is among the most used. According to \cite{cheng2009classifying}, the term K-means was suggested by \cite{macqueen1967some} to refer to algorithms that assign each item to the cluster with the nearest centroid (mean). Broadly speaking, this kind of algorithm usually starts with an initial partition into $K$ clusters given by the user and with the calculation of the centroids. Another option is to start with $K$ centroids prescribed by the user. Then, an item is re-assigned to the cluster having the nearest centroid and the centroid of the cluster that receives or loses the item is re-calculated. This procedure is repeated until no re-assignments are done. We remark that an important drawback of the K-means algorithm is that the clustering result might differ depending on the choice of the initial partition or initial centroids \cite{sauglam2006mixed}.

There are also optimization-based approaches to the clustering problem \cite{rao1971cluster}, which use mathematical programming models and involve exact or heuristic algorithms to solve the associated optimization problem. This type of approach is usually possible in small instances, due to the computational complexity of the related optimization problems \cite{sauglam2006mixed}. Dealing with this drawback is crucial in the case of the \textit{customer segmentation} problem because the most relevant real applications usually involve large datasets containing the data of thousands or even millions of customers, which could lead to problems that are computationally impossible to solve.

In this paper, we contribute by proposing an efficient optimization based approach in which solving the customer segmentation problem corresponds to clustering the vertices of a graph $G$, each of which representing a customer, in such a way that the vertices corresponding to customers with similar behaviors in terms of the $(R,F,M)$-score are in the same group. We represent this problem by a classical combinatorial optimization formulation known as the maximum-$k$-cut problem (or max-$k$-cut problem for short), where $k$ should be chosen as the number of groups.

To solve this max-$k$-cut problem in an acceptable computer processing time, even for large datasets of customers, we propose an auxiliary max-$k$-cut problem that is a reduction of the original one and whose associated graph $G^{\prime}$ has at most $T^{3}$ vertices, independently of the number of customers. Optimally solving the auxiliary problem leads us to the optimal solution to the original one as well, in such a way that we guarantee that customers with the same $(R,F,M)$-score are assigned to the same group. In the forthcoming Theorem \ref{teo_reduced_problem}, we prove this correspondence between these problems. The mathematical model of the max-$k$-cut problem for the graph $G^{\prime}$ is solved using Gurobi software. Our method has been applied to the Online Retail II dataset from the UCI Machine Learning Repository to corroborate our theoretical results and its effectiveness. Moreover, to the best of our knowledge, our proposal connecting graph theory, $(R,F,M)$-score and the max-$k$-cut problem to solve the customer segmentation is innovative and new in the literature.

The paper is organized as follows. In Section \ref{section:a graph based approach}, we present our graph-based approach to the customer segmentation problem, introducing graph $G$ and the mathematical formulation of the max-$k$-cut problem to be solved. We show how the reduced graph $G'$ is built from the original graph $G$ and how we use it to solve the max-$k$-cut problem. We also state and prove theoretical results that form the basis of our method. In Section \ref{sec: methodology}, we describe our graph-based methodology. In Section \ref{sec:results}, we perform computational experiments to show how our method works and then analyze the obtained results. Finally, Section \ref{sec: conclusion} presents our conclusions and perspective on future work.

%%%%%%%%%%%%%%%%%%%%%%%%%%%%%%%%%%%%%%%%%%%%%%%%%%%%%%%%%%%%%%%%%%%%%%%%%%%%%%%%%%%%%
%%%%%%%%%%%%%%%%%%%%%%%%%%%%%%%%%%%%%%%%%%%%%%%%%%%%%%%%%%%%%%%%%%%%%%%%%%%%%%%%%%%%%
%%%%%%%%%%%%%%%%%%%%%%%%%%%%%%%%%%%%%%%%%%%%%%%%%%%%%%%%%%%%%%%%%%%%%%%%%%%%%%%%%%%%%
%%%%%%%%%%%%%%%%%%%%%%%%%%%%%%%%%%%%%%%%%%%%%%%%%%%%%%%%%%%%%%%%%%%%%%%%%%%%%%%%%%%%%
%%%%%%%%%%%%%%%%%%%%%%%%%%%%%%%%%%%%%%%%%%%%%%%%%%%%%%%%%%%%%%%%%%%%%%%%%%%%%%%%%%%%%

\section{Customer segmentation: a graph-based approach}
\label{section:a graph based approach}

In this section, based on the $(R,F,M)$-scores of each customer, we show how customers are compared to each other and divided into groups using a graph approach. First, we need to introduce the max-$k$-cut problem.

\

%%%%%%%%%%%%%%%%%%%%%%%%%%%%%%%%%%%%%%%%%%%%%%%%%%%%%%%%%%%%%%%%%%%%%%%%%%%%%%%%%%%%%
%%%%%%%%%%%%%%%%%%%%%%%%%%%%%%%%%%%%%%%%%%%%%%%%%%%%%%%%%%%%%%%%%%%%%%%%%%%%%%%%%%%%%
%%%%%%%%%%%%%%%%%%%%%%%%%%%%%%%%%%%%%%%%%%%%%%%%%%%%%%%%%%%%%%%%%%%%%%%%%%%%%%%%%%%%%
%%%%%%%%%%%%%%%%%%%%%%%%%%%%%%%%%%%%%%%%%%%%%%%%%%%%%%%%%%%%%%%%%%%%%%%%%%%%%%%%%%%%%
%%%%%%%%%%%%%%%%%%%%%%%%%%%%%%%%%%%%%%%%%%%%%%%%%%%%%%%%%%%%%%%%%%%%%%%%%%%%%%%%%%%%%

\subsection{The max-\texorpdfstring{$k$}{k}-cut problem}

Consider a simple weighted graph denoted by $G=(V,E)$, where $V = \{ 1, 2, \dots, n \}$ is the set of vertices and $E$ is the set of edges $e_{ij}$ connecting vertices $i, j \in V$, with $|E| = m$. Each edge $e_{ij} \in E$ has a non-negative weight $w_{ij}$. Given the parameter $k \geq 2$, the max-$k$-cut problem consists of partitioning the set of vertices $V$ into $k$ disjoint subsets, say $V_1, V_2, \ldots, V_k$, so that the sum of the weights of the edges connecting vertices belonging to different subsets of the partition is the largest possible. It is an NP-hard combinatorial optimization problem \cite{van2016new}. 

Variations of this problem can be obtained by imposing constraints on each part $V_i$, for $i \in \{ 1, 2, \dots, k \}$. For instance: (i) if each $V_i$ induces a connected subgraph, then we have the connected max-$k$-cut problem \cite{healy2024branch,hojny2021mixed}; (ii) if $|V_i|$ is at most a given integer number, then we are dealing with the capacitated max-$k$-cut problem \cite{gaur2008capacitated}. Also, a wide variety of practical applications can be tackled via the max-$k$-cut problem and its variations, such as in the electricity market \cite{ambrosius2020endogenous,hojny2021mixed}, forest planning \cite{carvajal2013imposing}, scheduling \cite{carlson1966scheduling}, and clustering \cite{poland2006clustering}. Concerning the strategies to solve the max-$k$-cut problem, we refer the reader to the recent papers \cite{fakhimi2023relaxations,healy2024branch,hojny2021mixed} for approaches based on mathematical formulations and relaxations. Heuristic methods have also been proposed to provide good solutions to the max-$k$-cut problem \cite{gaur2008capacitated,ma2017multiple,zhu2013max}.

We consider the following mathematical formulation of the max-$k$-cut problem in terms of a binary quadratic optimization (BQO) problem. The decision variables $x_{il}$ represent whether a vertex $i$ belongs to or not to $V_l$, that is, $x_{il} = 1$ if vertex $i \in V_{l}$ and $x_{il} = 0$ otherwise. The mathematical  formulation \eqref{bqo_formulation} is

\begin{equation} \tag{BQO} \label{bqo_formulation}
  \begin{array}{rl}
    \max & \, f(x) = \D{\sum_{e_{ij} \in E}{w_{ij}} \left ( 1 - \sum_{l = 1}^k x_{il}x_{jl} \right )} \\[12pt]
    s.t. & \, \D{\sum_{l=1}^k x_{il}} = 1, \, i = 1,2, \dots, n \\[12pt]
    & \, x_{il} \in \{ 0, 1 \}, \mbox{ for } i=1,2,\dots,n \mbox{ and } l = 1, 2, \dots, k.
  \end{array}
\end{equation}
The maximized objective function considers only the weights related to the edges that connect vertices belonging to different subsets of the partition. The first set of constraints imposes that each vertex of the graph should be assigned to exactly one group, while the second set imposes that every $x_{il}$ is a binary variable.

The input for the mathematical model \eqref{bqo_formulation} is a weighted graph $G$ represented by its weighted adjacency matrix. Since we are interested in making a customer segmentation by using the RFM analysis, this graph is built as follows: each vertex $i$ of the graph $G$ represents a customer, which has a three-dimensional score $(R_i, F_i, M_i)$; two vertices $i$ and $j$ are connected by an edge of weight $w_{ij} \geq 0$ given by the Manhattan distance between the respective scores $(R_i, F_i, M_i)$ and $(R_j, F_j, M_j)$, that is,
\begin{equation*}
  w_{ij} = |R_i - R_j| + |F_i - F_j| +|M_i - M_j|.
\end{equation*}
Customers $i$ and $j$ with the same $(R,F,M)$-scores have $w_{ij} = 0$, which means that an edge does not connect the corresponding vertices.

\
%%%%%%%%%%%%%%%%%%%%%%%%%%%%%%%%%%%%%%%%%%%%%%%%%%%%%%%%%%%%%%%%%%%%%%%%%%%%%%%%%%%%%
%%%%%%%%%%%%%%%%%%%%%%%%%%%%%%%%%%%%%%%%%%%%%%%%%%%%%%%%%%%%%%%%%%%%%%%%%%%%%%%%%%%%%
%%%%%%%%%%%%%%%%%%%%%%%%%%%%%%%%%%%%%%%%%%%%%%%%%%%%%%%%%%%%%%%%%%%%%%%%%%%%%%%%%%%%%
%%%%%%%%%%%%%%%%%%%%%%%%%%%%%%%%%%%%%%%%%%%%%%%%%%%%%%%%%%%%%%%%%%%%%%%%%%%%%%%%%%%%%
%%%%%%%%%%%%%%%%%%%%%%%%%%%%%%%%%%%%%%%%%%%%%%%%%%%%%%%%%%%%%%%%%%%%%%%%%%%%%%%%%%%%%

\subsection{Reduction of graph \texorpdfstring{$G$}{G}: an auxiliary problem}
\label{subsec_graph_G_prime}

The number of vertices of $G$ is equal to the number of customers, and since a dataset can contain thousands of customers, finding the optimal solution to the associated max-$k$-cut problem (with a high number of vertices and edges) using model \eqref{bqo_formulation} may be computationally infeasible. For this reason, we propose constructing another graph from $G$, denoted here by $G^{\prime}$, which has fewer vertices and edges than $G$. The  construction of $G^{\prime}$ is as follows.

\

\begin{proc}{\bf (Reducing graph $G$)} \label{procedure_1}
\begin{itemize}
    \item \textit{ Merging vertices: Let $\{i_1, i_2,\ldots,i_p\}$ be the set of all vertices with a certain $(R,F,M)$-score in $G$. Merge those vertices into one vertex $i$ and add $i$ to the vertex set of $G^{\prime}$, denoted by $V(G^{\prime})$. Apply this procedure to all sets of vertices with the same $(R,F,M)$-scores in $G$. At the end, assume that $G^{\prime}$ has $n^{\prime}$ vertices (note that $n' \leq n$).}
    
    \item \textit{Updating weights: Let $i \in V(G^{\prime})$ obtained as a merge of $\{i_1, i_2,\ldots,i_p\}$, and let $j \in V(G^{\prime})$ obtained as a merge of $\{j_1, j_2,\ldots,j_q\}$. Then the edge weight connecting $i$ to a vertex $j$ in $G^{\prime}$ is given by $$w'_{ij} = \sum_{k=1}^{p} \sum_{r = 1}^{q} w_{i_k \, j_r}.$$ Apply this procedure to all vertices of $G^{\prime}.$ }
\end{itemize}
\end{proc}

\

\begin{figure}[H]
    \includegraphics[width=14cm]{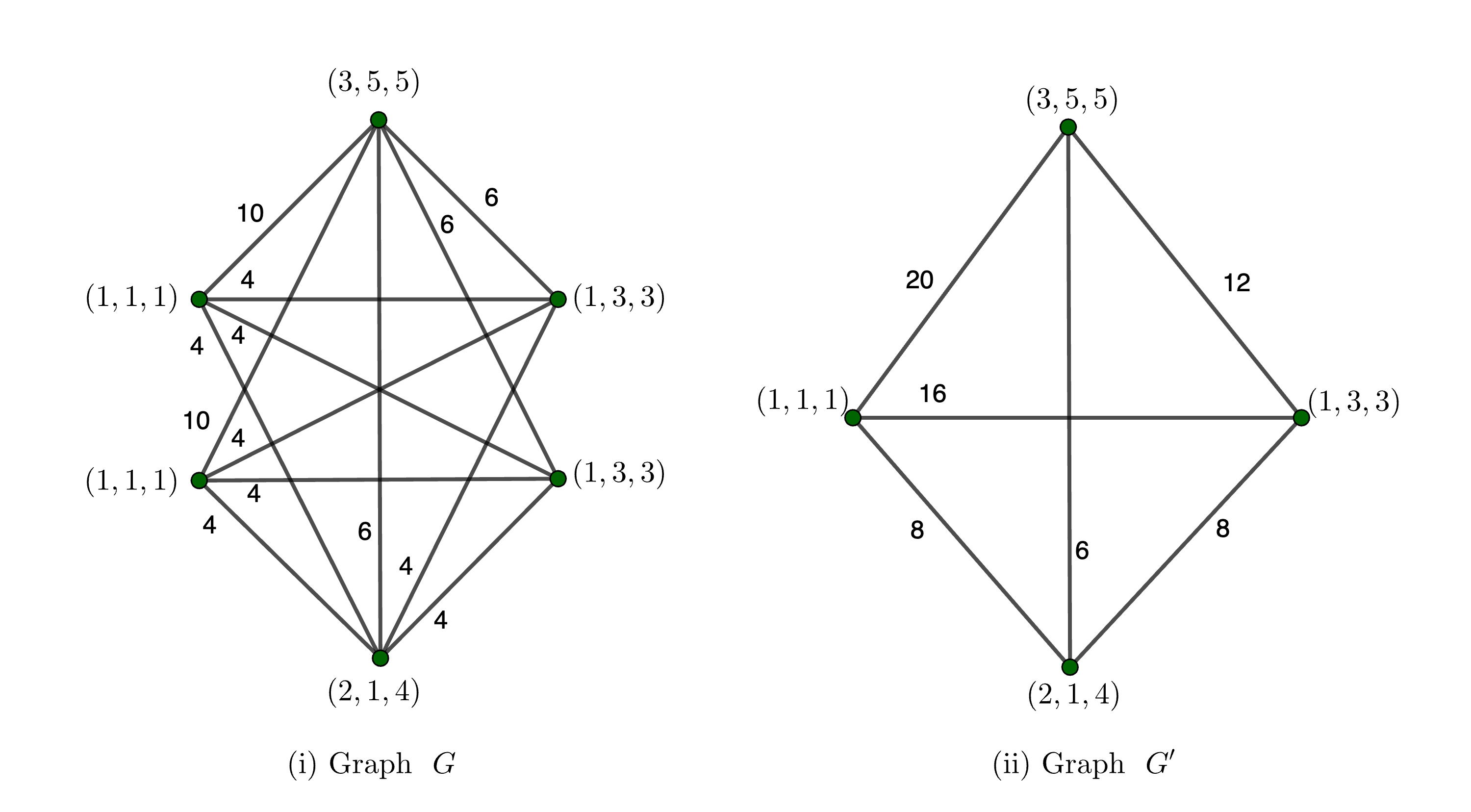}
    \centering
    \caption{Example of a graph $G^{\prime}$ in (ii) obtained from $G$ in (i).}
    \label{fig:reduction}
  \end{figure}

Adopting $T = 5$, consider graph $G$ in Figure \ref{fig:reduction} with 6 vertices and 26 edges. The vertices $(1,1,1)$ have no edge between them because $w_{(1,1,1),(1,1,1)} = 0.$ The graph $G^{\prime}$ is obtained by merging the vertices with the same $(R,F,M)$-scores $(1,1,1)$ and $(1,3,3)$. In the sequel, the weights of the corresponding incident edges are updated.

Note that, independently of the number of customers in the dataset, we ensure that graph $G^{\prime}$ has at most $125$ vertices, since there exist at most $5^3 = 125$ different $(R,F,M)$-scores. It is important to remark that, as we are going to prove, from a solution to the max-$k$-cut problem related to the graph $G^{\prime}$, we can easily obtain a solution to the problem related to $G$, but with the advantage that the mathematical model \eqref{bqo_formulation}, when obtained from $G^{\prime}$, may have considerably fewer variables and constraints and should be easier to solve than the original problem related to $G$. 

\
%%%%%%%%%%%%%%%%%%%%%%%%%%%%%%%%%%%%%%%%%%%%%%%%%%%%%%%%%%%%%%%%%%%%%%%%%%%%%%%%%%%%%
%%%%%%%%%%%%%%%%%%%%%%%%%%%%%%%%%%%%%%%%%%%%%%%%%%%%%%%%%%%%%%%%%%%%%%%%%%%%%%%%%%%%%
%%%%%%%%%%%%%%%%%%%%%%%%%%%%%%%%%%%%%%%%%%%%%%%%%%%%%%%%%%%%%%%%%%%%%%%%%%%%%%%%%%%%%
%%%%%%%%%%%%%%%%%%%%%%%%%%%%%%%%%%%%%%%%%%%%%%%%%%%%%%%%%%%%%%%%%%%%%%%%%%%%%%%%%%%%%
%%%%%%%%%%%%%%%%%%%%%%%%%%%%%%%%%%%%%%%%%%%%%%%%%%%%%%%%%%%%%%%%%%%%%%%%%%%%%%%%%%%%%

\subsection{Retrieving an optimal solution to the original problem}
\label{subsec_retrieve_solution}

Once we have reduced the graph $G$ to the graph $G^{\prime}$ employing Procedure \ref{procedure_1}, we solve problem \eqref{bqo_formulation} with graph $G^{\prime}$ as input, obtaining an optimal solution $\overline{x}^* \in \mathbb{R}^{n' \times k}$, where $n'$ is the number of vertices of $G^{\prime}$. Then, from $\overline{x}^*$, we can generate an optimal solution $x^* \in \mathbb{R}^{n \times k}$ to the problem \eqref{bqo_formulation} with the graph $G$ as input through the following procedure.

\vspace{0.5cm}

\begin{proc} \label{procedure_2}
  \textit{For each vertex $p$ in the graph $G^{\prime}$, if the vertex $i$ of the graph $G$ has been merged to the vertex $p$ in $G^{\prime}$ by Procedure \ref{procedure_1}, put $x^*_{il} = \overline{x}^*_{pl}$, for every $l \in \{1, \ldots, k \}$.}
\end{proc}

Additional comments about the relation between $G$ and $G^{\prime}$ for the customer segmentation problem are worth mentioning. When $G$ is given as input to the max-$k$-cut problem \eqref{bqo_formulation}, we expect that vertices with equal $(R,F,M)$-scores are assigned to the same group (or cluster) in the optimal solution. Unfortunately, it is not generally possible to guarantee that an optimal solution to the max-$k$-cut for the graph $G$ has this property. However, we prove that such a solution always exists and, that our technique always leads to this kind of solution, for which we give a more precise definition in the sequel to simplify the presentation of our arguments throughout the paper.

\begin{defi}{\bf (Segmentation property)} \label{defi_segmentation-property}
We say that a partition of the vertices of $G$ has the segmentation property if the vertices having the same $(R,F,M)$-score lie in the same group.
\end{defi}

\begin{defi}{\bf (Segmentation solution)} \label{defi_segmentation-solution}
A segmentation solution to problem \eqref{bqo_formulation} is a feasible solution that has the segmentation property.
\end{defi}

\

The existence of an optimal segmentation solution is proved in Theorem \ref{teo_segmentation_solution} by showing that, if you find an optimal solution violating this property, it is possible to make some local changes in this solution in such a way that we keep the optimal objective function value and rearrange vertices with equal $(R,F,M)$-scores to the same group, i.e., we obtain a new optimal solution that respects the fact that vertices with equal $(R,F,M)$-scores are always in the same group.

\begin{teo} \label{teo_segmentation_solution}
   Consider the graph $G$, built as described in Section 2.1, given as input to \eqref{bqo_formulation}. Then, there is at least one optimal segmentation solution for \eqref{bqo_formulation}.
\end{teo}
\begin{proof}[\bf Proof]
  Suppose that $x^*$ is an optimal solution to \eqref{bqo_formulation}. If $x^*$ is a segmentation solution, then the result is proven. Now, consider that $x^*$ is not a segmentation solution, that is, there are vertices with the same $(R,F,M)$-score but belonging to different partition sets.
  
  Let $(\overline{R},\overline{F},\overline{M})$ be an arbitrary score for which this happens, and let us fix a set $V_p$ of the partition containing a vertex associated with the score $(\overline{R},\overline{F},\overline{M})$. We denote this vertex associated with the score $(\overline{R},\overline{F},\overline{M})$, that lies in the set $V_p$, by $i$. We will conclude that any other vertex associated with the score $(\overline{R},\overline{F},\overline{M})$ can be moved to the set $V_p$ and we still have a global optimal solution. This way, we can find a global optimal solution in which all the vertices associated with the the score $(\overline{R},\overline{F},\overline{M})$ are grouped in the set $V_p$.
  
  Indeed, let $j$ be any other arbitrary vertex associated with the score $(\overline{R},\overline{F},\overline{M})$ that lies in a set $V_q$ of the partition, distinct from the set $V_p$. Let $S$ be the sum of the weights of the edges connecting vertices that lie in different sets of the partition but are not related to $i$ or $j$. Let $S_p$ be the sum of the weights of the edges connecting the vertex $i$ with vertices that lie in a set distinct from the set $V_p$ of the partition and let $S_q$ be the sum of the weights of the edges connecting the vertex $j$ with vertices that lie in a set distinct from the set $V_q$ of the partition. Then the value of the objective function in this optimal solution is $S + S_p + S_q$. 
  
  Considering that $i$ and $j$ are not connected by an edge and share the same neighbors, if we move the vertex $i$ from the set $V_p$ to the set $V_q$, the value of the objective function in the resulting partition is $S + S_q + S_q$. Since $S + S_p + S_q$ is the value of the objective function in a global optimal solution, we conclude that
  \begin{equation*}
    S + S_q + S_q \leq S + S_p + S_q,
  \end{equation*}
  which implies that $S_q \leq S_p$. Analogously, we can move the vertex $j$ from the set $V_q$ to the set $V_p$. In this case, the value of the objective function in the resulting partition is $S + S_p + S_p$, leading us to $S_p \leq S_q$. Therefore, we conclude that $S_p = S_q$.
  
  Hence, in particular, if we move the vertex $j$ to the set $V_p$ of the partition, the resulting partition, now with $i$ and $j$ grouped in the set $V_p$, is a global optimal solution. Since the vertex $j$ is arbitrary, any other vertex associated with the score $(\overline{R},\overline{F},\overline{M})$ can be grouped in the set $V_p$ and, at the end of this process, we obtain a global optimal solution with all the vertices associated with the score $(\overline{R},\overline{F},\overline{M})$ grouped in the same set.
  
  Considering now the resulting global optimal solution and observing that the score $(\overline{R},\overline{F},\overline{M})$ is also arbitrary, we can repeat this process with any other $(R,F,M)$-score for which we have associated vertices lying in different sets of the partition. This way, at the end of this process, we finally obtain a global optimal solution such that, for each fixed score $(R,F,M)$, all the vertices associated with this score $(R,F,M)$ are grouped in the same set, that is, we obtain a global optimal solution that is also a segmentation solution.
\end{proof}

\

Let us denote by $f$ the objective function of the \eqref{bqo_formulation} problem when $G$ is given as input, while we denote by $f_R$ the objective function of the \eqref{bqo_formulation} problem when $G^{\prime}$ is given as input. We have the following result.

\begin{teo} \label{teo_reduced_problem}
  Let ${x}^*$ be the optimal solution to the \eqref{bqo_formulation} problem when $G$ is given as input, and let $x_R^*$ be the optimal solution to the \eqref{bqo_formulation} problem when $G^{\prime}$ is given as input. Then, $f({x}^*) = f_R(x_R^*)$.
\end{teo}
\begin{proof}[\bf Proof] From Theorem \ref{teo_segmentation_solution}, we can assume that $x^{*}$ is a segmentation solution. It means that all vertices with the same $(R,F,M)$-scores are in the same group. Applying Procedure \ref{procedure_1} to graph $G$ with the corresponding solution $x^{*},$ we obtain a feasible solution $\hat{x}_R$ to the reduced problem, and by its construction, it is easy to see that $f(x^{*}) = f_R(\hat{x}_R)$.

Let $x^0$ be a feasible solution to the original problem obtained from $x_R^*$ by Procedure \ref{procedure_2}. By its construction, we have $f_R(x_R^*) = f(x^0)$. Accounting for the maximality of $f(x^*)$, we have $f(x^*) \geq f(x^0)$ and hence 
$$f_R(\hat{x}_R) = f(x^*) \geq f(x^0) = f_R(x_R^*).$$
By the maximality of $f_R(x_R^*)$, we then have
$$f_R(x_R^*) = f_R(\hat{x}_R),$$
from which we conclude that 
$f(x^*) = f(x^0) = f_R(x_R^*).$
\end{proof}

\

\begin{obsv}
  Note that, by the proof of Theorem \ref{teo_reduced_problem}, we also conclude that applying Procedure \ref{procedure_2} to an optimal solution $x_R^*$ of the reduced problem yields an optimal solution $x^0$ to the original problem.
\end{obsv}

%%%%%%%%%%%%%%%%%%%%%%%%%%%%%%%%%%%%%%%%%%%%%%%%%%%%%%%%%%%%%%%%%%%%%%%%%%%%%%%%%%%%%
%%%%%%%%%%%%%%%%%%%%%%%%%%%%%%%%%%%%%%%%%%%%%%%%%%%%%%%%%%%%%%%%%%%%%%%%%%%%%%%%%%%%%
%%%%%%%%%%%%%%%%%%%%%%%%%%%%%%%%%%%%%%%%%%%%%%%%%%%%%%%%%%%%%%%%%%%%%%%%%%%%%%%%%%%%%
%%%%%%%%%%%%%%%%%%%%%%%%%%%%%%%%%%%%%%%%%%%%%%%%%%%%%%%%%%%%%%%%%%%%%%%%%%%%%%%%%%%%%
%%%%%%%%%%%%%%%%%%%%%%%%%%%%%%%%%%%%%%%%%%%%%%%%%%%%%%%%%%%%%%%%%%%%%%%%%%%%%%%%%%%%%

\section{Methodology}
\label{sec: methodology}

Consider a dataset composed of $n$ customers. We propose a method for \textit{customer segmentation} using the RFM analysis and clustering by solving the max-$k$-cut problem related to the dataset. In the sequel, we describe our graph-based methodology.

\begin{itemize}
  \item Split the data (recency, frequency, and monetary value), sorting them in ascending order into $T$ equally spaced intervals, and assign an integer in the range 1 to $T$ to each variable $R_{i}, F_i, M_i$ for each customer $i \in \{1,\ldots,n\}$.
  \item Determine the possible values of $k$ (the number of clusters in which the customers will be segmented) according to the appropriate criteria, depending on the particular application.
  \item Build the corresponding weighted graph $G$ on $n$ vertices with the customer dataset. Each customer $i$ is a vertex of $G$, and two customers $i$ and $j$ are connected by an edge with weight $w_{ij} = |R_i-R_j|+|F_i-F_j|+ |M_i-M_j|.$
  \item Reduce $G$ to the graph $G^{\prime}$ by merging vertices with the same $(R,F,M)$-score and updating the edge weights as previously described in Procedure \ref{procedure_1}. The graph $G^{\prime}$ has $n^{\prime} \leq n$ vertices.
  \item For each one of the selected values of $k$, solve the corresponding max-$k$-cut problem through the formulation \eqref{bqo_formulation}, given $G^{\prime}$ with $n'$ vertices as input. Let $\overline{x}^* \in \mathbb{R}^{n' \times k}$ be the solution found.
  \item Obtain a segmentation solution $x^* \in \mathbb{R}^{n \times k}$ to the problem \eqref{bqo_formulation} with graph $G$ given as input using Procedure \ref{procedure_2}.
  \item Choose the best value of $k$ through some measure of clustering performance \cite{karanikola2021investigating} or other criteria that may be more suitable for a given practical application.
\end{itemize}

\

%%%%%%%%%%%%%%%%%%%%%%%%%%%%%%%%%%%%%%%%%%%%%%%%%%%%%%%%%%%%%%%%%%%%%%%%%%%%%%%%%%%%%
%%%%%%%%%%%%%%%%%%%%%%%%%%%%%%%%%%%%%%%%%%%%%%%%%%%%%%%%%%%%%%%%%%%%%%%%%%%%%%%%%%%%%
%%%%%%%%%%%%%%%%%%%%%%%%%%%%%%%%%%%%%%%%%%%%%%%%%%%%%%%%%%%%%%%%%%%%%%%%%%%%%%%%%%%%%
%%%%%%%%%%%%%%%%%%%%%%%%%%%%%%%%%%%%%%%%%%%%%%%%%%%%%%%%%%%%%%%%%%%%%%%%%%%%%%%%%%%%%
%%%%%%%%%%%%%%%%%%%%%%%%%%%%%%%%%%%%%%%%%%%%%%%%%%%%%%%%%%%%%%%%%%%%%%%%%%%%%%%%%%%%%

\section{Results and discussion} \label{sec:results}

Our computational experiments were performed on an Intel Core i$7$ processor operating at $3.4$ GHz, equipped with $64$ GB of RAM and running the Ubuntu Linux operating system. The algorithms were coded in Python 3.12.2, and we used Gurobi Optimizer version $11.0.2$ as the \eqref{bqo_formulation} solver in its default settings.

\

%%%%%%%%%%%%%%%%%%%%%%%%%%%%%%%%%%%%%%%%%%%%%%%%%%%%%%%%%%%%%%%%%%%%%%%%%%%
%%%%%%%%%%%%%%%%%%%%%%%%%%%%%%%%%%%%%%%%%%%%%%%%%%%%%%%%%%%%%%%%%%%%%%%%%%%
%%%%%%%%%%%%%%%%%%%%%%%%%%%%%%%%%%%%%%%%%%%%%%%%%%%%%%%%%%%%%%%%%%%%%%%%%%%
%%%%%%%%%%%%%%%%%%%%%%%%%%%%%%%%%%%%%%%%%%%%%%%%%%%%%%%%%%%%%%%%%%%%%%%%%%%
%%%%%%%%%%%%%%%%%%%%%%%%%%%%%%%%%%%%%%%%%%%%%%%%%%%%%%%%%%%%%%%%%%%%%%%%%%%

\subsection{Computational experiments} \label{subsec_computational}

In these tests, we utilized the Online Retail II dataset from the UCI Machine Learning Repository. This widely recognized dataset, employed in several studies, including \cite{rahim2021rfm} and \cite{rungruang2024rfm}, encompasses transactions from a non-store online retail operation based in the UK, specialized in the sale of unique all-occasion giftware, primarily for wholesalers. The period considered spans from December $1$, $2009$ to December $9$, $2011$, totaling $1,067,371$ transactions and comprising eight variables, including Customer ID, InvoiceDate, Quantity, and Price.

Before applying our method, we preprocessed the data. During the preprocessing phase, we removed entries with missing Customer IDs, negative values in Quantity and/or Price, as well as duplicate records. As a result, the dataset was reduced to $791,045$ transactions made by $5,878$ unique customers.

After data preprocessing, the RFM analysis was performed for each of the $5,878$ customers, adopting $T = 5$. For each customer $i$, the frequency ($F_i$) is defined as the number of transactions made by $i$ in the considered period; the monetary value ($M_i$) is calculated as the sum of the monetary transaction values (quantity x price) for $i$; and the recency ($R_i$) represents the elapsed time (in days) between the date of the last purchase made by $i$ and the last purchase in the entire dataset. Recency was transformed so that a lower value indicates a better outcome, aligning it with the other variables (frequency and monetary value).

Thus, each customer was assigned an $(R, F, M)$-score, with the three variables (R, F, and M) receiving a score from $1$ to $5$ based on the quintiles (division of the dataset into five equal parts) of each variable. A score of $1$ represents the worst case, while $5$ represents the best case.

To illustrate the difference between solving the formulation \eqref{bqo_formulation} for the original problem, in which the graph $G$ is given as input, and applying our method, in which the reduced graph $G^{\prime}$ is given as input, we first solved three small instances of the formulation \eqref{bqo_formulation} to their optimality using both the original and the reduced problems. The instances contain, respectively, the first $40$, $45$, and $50$ customers of the dataset resulting from the aforementioned preprocessing step.
 
The results for each one of these instances are displayed in Tables \ref{tab:tab1}, \ref{tab:tab2}, and \ref{tab:tab3}, where the columns are: the number $k$ of groups, for $k = 2$, $3$, and $4$; the average execution time (in seconds), considering $10$ executions of the formulation \eqref{bqo_formulation} when $G$ is given as input ($t_{G}$) and when $G^{\prime}$ is given as input $(t_{G^{\prime}})$, respectively; and the optimal objective function value, which was the same for both the original and the reduced problems, as it has been previously proved in Theorem \ref{teo_reduced_problem}.

\begin{table}[H]
\centering
\begin{tabular}{| c | c | c | c |}
 \hline
 Number of groups ($k$) & $t_{G}$ $(s)$ & $t_{G^{\prime}}$ $(s)$ & Objective function \\[0.5ex] 
 \hline
 $2$ & $0.90$ & $0.20$ & $2,428$ \\[0.5ex]
 \hline
 $3$ & $88.51$ & $3.15$ & $2,864$ \\[0.5ex]
 \hline
 $4$ & $314.42$ & $61.71$ & $3,052$ \\[0.5ex]
 \hline
\end{tabular}
\caption{Clustering results for the instance with $40$ customers for $k = 2$, $3$ and $4$, when graphs $G$ ($40$ vertices and $764$ edges) and $G^{\prime}$ ($28$ vertices and $378$ edges) are given as input to the formulation \eqref{bqo_formulation}.}
\label{tab:tab1}
\end{table}

\begin{obsv}
  The optimal objective function value of the formulation \eqref{bqo_formulation} increases as the number $k$ of clusters increases, but it does not have any implications for the clustering quality. This happens because a partition of the vertices into $k$ distinct groups is a particular case of a partition of the vertices into $k + 1$ groups, one of them having no vertices assigned to it. Hence, when we look for the optimal objective function value with more clusters available, more vertices can be assigned to distinct clusters and thus more edge weights can be considered in the objective function of the formulation \eqref{bqo_formulation}, increasing its optimal value, when compared to the optimal value when a smaller value of $k$ is considered.
\end{obsv}

\

\begin{table}[H]
\centering
\begin{tabular}{| c | c | c | c |} 
 \hline
 Number of groups ($k$) & $t_G$ ($s$)& $t_{G^{\prime}}$ ($s$) & Objective function \\[0.5ex] 
 \hline
 $2$ & $0.04$ & $ 0.18$ & $3,132$ \\[0.5ex]
 \hline
 $3$ & $142.10$ & $6.27$ & $3,702$ \\[0.5ex]
 \hline
 $4$ & $1,559.78$ & $59.58$ & $3,932$ \\[0.5ex]
 \hline
\end{tabular}
 \caption{Clustering results for the instance with $45$ customers for $k = 2$, $3$, and $4$, when graphs $G$ ($45$ vertices and $972$ edges) and $G^{\prime}$ ($32$ vertices and $496$ edges) are given as input to the formulation \eqref{bqo_formulation}.}
\label{tab:tab2}
\end{table}

\begin{table}[H]
\centering
\begin{tabular}{| c | c | c | c |}
 \hline
 Number of groups ($k$) & $t_G$ (s) & $t_{G^{\prime}}$ (s) & Objective function \\[0.5ex] 
 \hline
 $2$ & $0.11 $ & $0.22 $ & $3,930$ \\[0.5ex]
 \hline
 $3$ & $510.55 $ & $3.12 $ & $4,646$ \\[0.5ex]
 \hline
 $4$ & $5,245.87 $ & $70.83 $ & $4,948$ \\[0.5ex]
 \hline
\end{tabular}
\caption{Clustering results for the instance with $50$ customers for $k = 2$, $3$, and $4$, when graphs $G$ ($50$ vertices and $1,201$ edges) and $G^{\prime}$ ($33$ vertices and $528$ edges) are given as input to the formulation \eqref{bqo_formulation}.}
\label{tab:tab3}
\end{table}

From these tests with the small instances, we can see that even for small problems, involving datasets with few customers, our strategy may present a great advantage when compared to solving the formulation \eqref{bqo_formulation} with the original graph $G$. Indeed, for values of $k$ greater than $2$, our method needed a much lower computer processing time to solve the formulation \eqref{bqo_formulation}. Moreover, our approach guarantees that the solution found is a segmentation solution.

In Table \ref{tab:tab4}, the results considering all $5,878$ customers are presented. In this case, the graph $G$, for the original problem, has $5,878$ vertices and $16,874,125$ edges, and it was not possible to run the formulation \eqref{bqo_formulation} on our computer, for any value of $k$. As expected for these big problems, the formulation \eqref{bqo_formulation} can be solved only using our approach, using the graph $G^{\prime}$, which corresponds to the reduced problem and has $114$ vertices and $6,441$ edges. The execution time limit was set to two hours, or $7,200$ seconds.

The first column of Table \ref{tab:tab4} brings the number $k$ of groups, from $k = 2$ to $10$. In the second column, we have the execution time (in seconds) of the formulation \eqref{bqo_formulation} when $G^{\prime}$ is given as input. The execution time of $7,200 \, s$ means that it stopped because it reached the time limit, not the optimality criteria, which was reached only in the case $k = 2$. The objective function values found are listed in the third column of Table \ref{tab:tab4}. Finally, the overall silhouette index for each value of $k$ is presented in the fourth column.

The silhouette index \cite{rousseeuw1987silhouettes} is a metric used to assess the quality of a clustering by considering both the internal cohesion of the clusters and the separation between them. For each vertex, it calculates how close this vertex is to others in the same cluster compared to vertices in neighboring clusters, giving us the individual silhouette index of each vertex. The average of these individual values yields the overall silhouette score. The closer this average value is to 1, the better the clustering, indicating that the data is well grouped and well separated. Values close to 0 suggest overlapping clusters, and negative values indicate poor clustering.

\begin{table}[H]
\centering
\begin{tabular}{| c | c | c | c |} 
 \hline
 Number of groups ($k$) & Time (s) &  Objective function & Silhouette index \\[0.5ex] 
 \hline
 $2$ & $73.74$ & $56,957,982$ & 0.5030 \\[0.5ex]
 \hline
 $3$ & $7,200$ & $68,685,769$ & 0.4050 \\[0.5ex]
 \hline
 $4$ & $7,200$ & $73,706,438$ & 0.4289 \\[0.5ex]
 \hline
 $5$ & $7,200$ & $76,353,219$ & 0.3960 \\[0.5ex]
 \hline
 $6$ & $7,200$ & $77,803,712$ & 0.3576 \\[0.5ex]
 \hline
 $7$ & $7,200$ & $78,790,442$ & 0.3625 \\[0.5ex]
 \hline
 $8$ & $7,200$ & $79,654,918$ & 0.3774 \\[0.5ex]
 \hline
 $9$ & $7,200$ & $80,190,479$ & 0.3865 \\[0.5ex]
 \hline
 $10$ & $7,200$ & $80,646,842$ & 0.4008 \\[0.5ex]
 \hline
\end{tabular}
\caption{Clustering results for the instance with $5,878$ customers for $k = 2$ to $10$, when graph $G^{\prime}$ ($114$ vertices and $6,441$ edges) is given as input to the formulation \eqref{bqo_formulation}.}
\label{tab:tab4}
\end{table}

As shown in Table \ref{tab:tab4}, the highest silhouette indexes are obtained when customers are segmented into $2$ and $4$ clusters, respectively, indicating that $k = 2$ and $k = 4$ yield the most effective clusterings. Although the highest silhouette index is reached for $k = 2$, as we mentioned at the end of Section \ref{sec: methodology}, additional criteria may be considered when evaluating clustering quality and selecting the number $k$ of clusters. For instance, in \cite{rungruang2024rfm}, the number of groups was determined based on whether the resulting clusters were meaningful and valuable from a business perspective.

\

%%%%%%%%%%%%%%%%%%%%%%%%%%%%%%%%%%%%%%%%%%%%%%%%%%%%%%%%%%%%%%%%%%%%%%%%%%%
%%%%%%%%%%%%%%%%%%%%%%%%%%%%%%%%%%%%%%%%%%%%%%%%%%%%%%%%%%%%%%%%%%%%%%%%%%%
%%%%%%%%%%%%%%%%%%%%%%%%%%%%%%%%%%%%%%%%%%%%%%%%%%%%%%%%%%%%%%%%%%%%%%%%%%%
%%%%%%%%%%%%%%%%%%%%%%%%%%%%%%%%%%%%%%%%%%%%%%%%%%%%%%%%%%%%%%%%%%%%%%%%%%%
%%%%%%%%%%%%%%%%%%%%%%%%%%%%%%%%%%%%%%%%%%%%%%%%%%%%%%%%%%%%%%%%%%%%%%%%%%%

\subsection{Business analysis} \label{subsec_business}

In our application, the forthcoming analysis suggests that segmenting customers into four distinct groups may provide more meaningful and actionable insights for business applications.

Table \ref{tab:tab5} presents a descriptive statistics analysis of $RFM$ for all $5,878$ customers.

\begin{table}[H]
\centering
\begin{tabular}{| c | c | c | c | c | c |}
 \hline
 Cluster & Value & Recency & Frequency & Monetary & Number of customers \\[0.5ex] 
 \hline
 \multirow{3}{*}{$1$} & min & $0.00$ & $1.00$ & $2.95$ & \multirow{3}{*}{$5,878$}\\[0.5ex]
 \cline{2-5}
 & mean & $200.87$ & $6.30$ & $2,965.56$ & \\[0.5ex]
 \cline{2-5}
 & max & $738.00$ & $398.00$ & $580,987.04$ & \\[0.5ex]
 \hline
\end{tabular}
\caption{Minimum, maximum and mean values of recency, frequency and monetary value for a unique group containing all $5,878$ customers.}
\label{tab:tab5}
\end{table}

%\obs{Depois para $k = 2$, explicando o resultado}

As expected, considering all the customers together, it is not clear how to gain insight into their individual behavior or their relationship with the company. Therefore, clustering is necessary in order to identify groups of customers with similar behaviors and then develop targeted marketing strategies for each group. In Table \ref{tab:tab6}, in the sequel, we show the results for $k = 2$, which yielded the highest silhouette index.

\begin{table}[H]
\centering
\begin{tabular}{| c | c | c | c | c | c |}
 \hline
 Cluster & Value & Recency & Frequency & Monetary & Number of customers \\[0.5ex] 
 \hline
 \multirow{3}{*}{$1$} & min & $0.00$ & $1.00$ & $2.95$ & \multirow{3}{*}{$2,898$}\\[0.5ex]
 \cline{2-5}
 & mean & $325.22$ & $1.66$ & $472.10$ & \\[0.5ex]
 \cline{2-5}
 & max & $738.00$ & $8.00$ & $13,916.34$ & \\[0.5ex]
 \hline
 \multirow{3}{*}{$2$} & min & $0.00$ & $1.00$ & $30.95$ & \multirow{3}{*}{$2,980$} \\[0.5ex]
 \cline{2-5}
 & mean & $79.93$ & $10.78$ & $5,390.40$ & \\[0.5ex]
 \cline{2-5}
 & max & $691.00$ & $398.00$ & $580,987.04$ & \\[0.5ex]
 \hline
\end{tabular}
\caption{Minimum, maximum and mean recency, frequency and monetary values for two clusters.}
\label{tab:tab6}
\end{table}
From Table \ref{tab:tab6}, we see that customers in Cluster $1$ take about $325$ days, on average, to make a new purchase, have low buying frequency (less than two purchases on average), and have low monetary value spent over the considered period. Therefore, they are currently customers who contribute a low income to the company.

On the other hand, on average, customers in Cluster $2$, when compared to customers in Cluster $1$, take considerably less time to make a new purchase, have higher buying frequency, and greater monetary value spent over the considered period. Therefore, they are currently a group of more loyal customers and generate higher income for the company.

However, even after identifying two very distinct groups of customers, some questions about the customers' behavior still remain. For instance, considering Cluster $1$, we are not able to differentiate customers who could be stimulated to make more frequent purchases from those who have possibly been lost to competitors. In Cluster $2$, we are not able to distinguish the most loyal and profitable customers from those who could be the target of marketing strategies to increase their buying frequency and the amount of money spent in the company.

We understand that these questions are better answered when we divide the customers into four distinct groups. Indeed, in this case, we are able to identify more customer segments, which helps us analyze their behavior and relationship with the company in more detail. The results for $k = 4$, which yielded the second highest silhouette index, are shown in Table \ref{tab:tab7}.

\begin{table}[H]
\centering
\begin{tabular}{| c | c | c | c | c | c |}
 \hline
 Cluster & Value & Recency & Frequency & Monetary & Number of customers \\[0.5ex] 
 \hline
 \multirow{3}{*}{$1$} & min & $58.00$ & $1.00$ & $2.95$ & \multirow{3}{*}{$1,604$}\\[0.5ex]
 \cline{2-5}
 & mean & $417.27$ & $1.19$ & $283.07$ & \\[0.5ex]
 \cline{2-5}
 & max & $738.00$ & $8.00$ & $2,803.20$ & \\[0.5ex]
 \hline
 \multirow{3}{*}{$2$} & min & $58.00$ & $1.00$ & $86.15$ & \multirow{3}{*}{$1,403$} \\[0.5ex]
 \cline{2-5}
 & mean & $282.17$ & $4.05$ & $1,493.49$ & \\[0.5ex]
 \cline{2-5}
 & max & $738.00$ & $17.00$ & $44,534.30$ & \\[0.5ex]
 \hline
 \multirow{3}{*}{$3$} & min & $0.00$ & $3.00$ & $808.62$ & \multirow{3}{*}{$1,612$}\\[0.5ex]
 \cline{2-5}
 & mean & $39.66$ & $16.13$ & $8,575.24$ & \\[0.5ex]
 \cline{2-5}
 & max & $576.00$ & $398.00$ & $580,987.04$ & \\[0.5ex]
 \hline
 \multirow{3}{*}{$4$} & min & $0.00$ & $1.00$ & $20.80$ & \multirow{3}{*}{$1,259$}\\[0.5ex]
 \cline{2-5}
 & mean & $40.95$ & $2.66$ & $841.01$ & \\[0.5ex]
 \cline{2-5}
 & max & $186.00$ & $15.00$ & $168,472.50$ & \\[0.5ex]
 \hline
\end{tabular}
\caption{Minimum, maximum and mean recency, frequency and monetary values for four clusters.}
\label{tab:tab7}
\end{table}
From Table \ref{tab:tab7}, we note that, in Cluster $1$, we have customers who have practically abandoned the company, because they have very high recency, very low frequency (on average one purchase), and low value spent. This group could be the target of campaigns to encourage their return, so the company can reactivate these customers.

In Cluster $2$, customers are slightly more frequent and spend more, but their last purchase ocurred some time ago. This group could be the target of retention campaigns or loyalty programs.

Cluster $3$, on the other hand, is the company's most valuable customer group. These customers make purchases almost every month, with high frequency and high value spent on transactions. Campaigns can target this group to maintain and eventually improve this relationship with the company.

Finally, in Cluster $4$, customers have low recency, indicating that they have bought recently, but with low frequency and low value spent. The company can develop marketing strategies aimed at this group so that they buy more frequently and increase the amount spent.

Therefore, by considering clustering quality metrics and the significant business implications, a company can choose the ideal number $k$ of clusters. Hence, by applying an effective technique to segment customers, such as the one presented in this research, a company leverage its business by no longer treating all customers equally, but instead considering them according to their purchasing profiles.

%%%%%%%%%%%%%%%%%%%%%%%%%%%%%%%%%%%%%%%%%%%%%%%%%%%%%%%%%%%%%%%%%%%%%%%%%%%%%%%%%%%%%
%%%%%%%%%%%%%%%%%%%%%%%%%%%%%%%%%%%%%%%%%%%%%%%%%%%%%%%%%%%%%%%%%%%%%%%%%%%%%%%%%%%%%
%%%%%%%%%%%%%%%%%%%%%%%%%%%%%%%%%%%%%%%%%%%%%%%%%%%%%%%%%%%%%%%%%%%%%%%%%%%%%%%%%%%%%
%%%%%%%%%%%%%%%%%%%%%%%%%%%%%%%%%%%%%%%%%%%%%%%%%%%%%%%%%%%%%%%%%%%%%%%%%%%%%%%%%%%%%
%%%%%%%%%%%%%%%%%%%%%%%%%%%%%%%%%%%%%%%%%%%%%%%%%%%%%%%%%%%%%%%%%%%%%%%%%%%%%%%%%%%%%

\section{Conclusion} \label{sec: conclusion}

%In the present article, we proposed a graph-based clustering method aimed at solving customer segmentation problems. Our approach combines the RFM analysis, a widely used \textit{customer value analysis} method which identifies the characteristics of customers using only the variables recency, frequency and monetary value, with the max-$k$-cut problem, a classical optimization problem which consists in partitioning the vertices of a weighted simple graph into $k$ disjoint sets so that a certain objective function is maximized in formulation \eqref{bqo_formulation}.

This article proposed a graph-based clustering method to solve customer segmentation problems. Our approach combines the RFM analysis with the max-$k$-cut problem. Each customer is represented by a vertex of a graph $G$, and the edge weights are the Manhattan distances between the vectors corresponding to the customers' $(R,F,M)$-scores. To segment the customers into $k$ disjoint groups, we must solve the corresponding formulation \eqref{bqo_formulation} of the max-$k$-cut problem.

%To apply our method to the customer segmentation problem, we first build the weighted graph $G$, in which each customer is represented by a vertex, and the edge weights are the Manhattan distances between the vectors corresponding to the customers' $(R,F,M)$-scores. To segment the customers into $k$ disjoint groups we then have to solve the corresponding max-$k$-cut problem. 

However, customer datasets may contain thousands of customers, and solving the formulation \eqref{bqo_formulation} when $G$ is given as input may be computationally infeasible. In fact, in our computational experiments, in Subsection \ref{subsec_computational}, the graph $G$ had $5,878$ vertices and $16,874,125$ edges, and it was not possible to run the corresponding max-$k$-cut problem on our computer, for any value of $k$.

This is why, in our method, we introduced the reduced graph $G^{\prime}$ of $G$, which has a considerably smaller number of vertices and edges, independent of the number of customers. We then analyzed the auxiliary problem obtained when $G^{\prime}$ is given as input in formulation \eqref{bqo_formulation}, proving, in Theorem \ref{teo_reduced_problem}, that, from any optimal solution to this problem, we can easily retrieve a segmentation solution to the original problem, when $G$ is given as input, via Procedure \ref{procedure_2}.

Next, to validate our method, we applied it to the Online Retail II dataset from the UCI Machine Learning Repository in Section \ref{sec:results}. After preprocessing the data, we obtained a dataset containing the $(R,F,M)$-scores of $5,878$ customers. The initial computational experiments with small instances, using only the first $40$, $45$, and $50$ customers of this dataset, already suggested the advantage of using our method instead of solving formulation \eqref{bqo_formulation} with the original graph given as input. Indeed, besides our method provides a segmentation solution for each problem, for $k > 2$, the computer processing time is significantly smaller when we use the reduced graph.

The full problem, with all $5,878$ customers, in its turn, was solvable only through our method, passing the reduced graph $G^{\prime}$ as input to the formulation \eqref{bqo_formulation}. We solved the associated max-$k$-cut problem for $k = 2$ to $10$ and calculated the silhouette index of the resulting clusterings of the $5,878$ customers, identifying that the highest silhouette indexes were obtained for $k = 2$ and $4$, in this order.

Nevertheless, to choose the number $k$ of clusters, we argue that other criteria, in addition to cluster quality metrics, should be considered, depending on the practical application we are working with. In fact, in Subsection \ref{subsec_business}, we analyzed the resulting clustering of the customers for $k = 2$ and $k = 4$ through the minimum, mean and maximum values of the variables recency, frequency and monetary value of each cluster.

The conclusion was that dividing the customers into four groups enabled us to identify other segments of interest that were not observable in the division into two groups and could be targeted by more personalized and effective marketing strategies. Therefore, in our case, despite $k = 2$ yielding the highest silhouette index, $k = 4$ should be preferred because the resulting clustering is more meaningful and valuable from the business perspective.

We highlight that the methodology presented in this article can be applied in a more general scenario. Suppose that one wishes to classify some items (customers, products, students, etc.) and, to do so, you have a system consisting of $q$ variables $Var_1$, $Var_2$, $\ldots$, $Var_q$, instead of $R$, $F$ and $M$, which can assume only integer values from $1$ to some fixed integer $T$. Our approach can be adapted entirely to this new situation, with the graph associated with the reduced problem having at most $T^q$ vertices, and our method should perform well in applications for which $T^q$ is not too large.

Examples of these more general scenarios are the models for customer segmentation with scales for the RFM variables adopting different values for $T$ \cite{dougan2018customer,panuvs2016customer,ramkumar2025enhancing} or which incorporate more variables to the RFM model, such as the LRFM, RFMTC, RFMD and RFMT models \cite{chang2004integrating,yeh2009knowledge,noori2015analysis,zhou2021customer}. Future works could focus on using our method in practical applications that employ these more general models.

Moreover, we highlight that, in this article, we tried to solve the max-$k$-cut problem to optimality using the Gurobi software. Therefore, future works can address the use of other methods to solve the max-$k$-cut problem, such as metaheuristics. These non-exact methods may be especially relevant if we consider models with more variables or scales with $T > 5$, which result in reduced graphs $G^{\prime}$ that have more vertices than we have in our case study, with three variables ($R$, $F$ and $M$) and $T = 5$.

\
%%%%%%%%%%%%%%%%%%%%%%%%%%%%%%%%%%%%%%%%%%%%%%%%%%%%%%%%%%%%%%%%%%%%%%%%%%%%%%%%%%%%%
%%%%%%%%%%%%%%%%%%%%%%%%%%%%%%%%%%%%%%%%%%%%%%%%%%%%%%%%%%%%%%%%%%%%%%%%%%%%%%%%%%%%%
%%%%%%%%%%%%%%%%%%%%%%%%%%%%%%%%%%%%%%%%%%%%%%%%%%%%%%%%%%%%%%%%%%%%%%%%%%%%%%%%%%%%%
%%%%%%%%%%%%%%%%%%%%%%%%%%%%%%%%%%%%%%%%%%%%%%%%%%%%%%%%%%%%%%%%%%%%%%%%%%%%%%%%%%%%%
%%%%%%%%%%%%%%%%%%%%%%%%%%%%%%%%%%%%%%%%%%%%%%%%%%%%%%%%%%%%%%%%%%%%%%%%%%%%%%%%%%%%%

\section{Data availability}

The Online Retail II dataset from the UCI Machine Learning Repository, which we used in our computational experiments, is available at \url{https://archive.ics.uci.edu}.

\vspace{1cm}

\noindent{\bf Acknowledgments.} The first author  would like to thank CAPES for the
support received during 2024/2025. The second author is supported by CNPq grants 315739/2021-5 and 403963/2021-4.\\

%% Use \subsubsection, \paragraph, \subparagraph commands to 
%% start 3rd, 4th and 5th level sections.
%% Refer following link for more details.
%% https://en.wikibooks.org/wiki/LaTeX/Document_Structure#Sectioning_commands

%\appendix
%\section{Maybe the proofs of the results as appendices? See how it is done in other articles of the Journal.}
%\label{app1}


%%       \cite{<label>} ==> [1]

%%

%% If you have bib dataset file and want bibtex to generate the
%%  \bibliography{<your bibdataset>}



\end{document}